\documentclass{conm-p-l}
\usepackage{amsmath}
\usepackage{ifthen}
\usepackage{amsfonts}
\usepackage{amssymb}
\usepackage{epsf}
\usepackage{amsopn}

\newcommand{\showcomments}{yes}
\newcommand{\refe}[1]{\eqref{e:#1}}

\newcommand{\comment}[1]
{\ifthenelse{\equal{\showcomments}{yes}}
{\footnotemark\marginpar{\sffamily{\tiny
\addtocounter{footnote}{-1}\footnotemark#1

}\normalfont}}{}}

\newtheorem{thm}{Theorem}[section]
\newtheorem{lem}[thm]{Lemma}

\newtheorem*{torsion rf theorem}{Theorem~\ref{thm:torsion rf}}

\theoremstyle{definition}
\newtheorem{defn}[thm]{Definition}

\newtheorem{notation}[thm]{Notation}

\newtheorem{quest}[thm]{Question}

\begin{document}

\title{On definitions of relatively hyperbolic groups}
\author{Inna Bumagin}
\address{Department of Mathematics and Statistics, McGill
University, Montr\'{e}al, Qu\'{e}bec H3A 2K6, Canada}

\email{bumagin@math.mcgill.ca}



\begin{abstract} The purpose of this note is to provide a short
alternate proof
of the fact that \cite[Question~1]{Sz} has an affirmative answer.
Our proof combined with the result of Szczepanski \cite{Sz} shows
that a group which is relatively hyperbolic in the sense of the
definition of Gromov is relatively hyperbolic in the sense of the
definition of Farb.
\end{abstract}

\maketitle

\section{Definitions}

The notion of a relatively hyperbolic group was introduced by
Gromov~\cite{Gr} as a generalization of the concept of a (word)
hyperbolic group.

\subsection{The Gromov definition}\label{defn:Gromov} Let $X$ be a hyperbolic
(in the Gromov sense) complete locally compact geodesic metric
space. Let $x\in X$, suppose $z$ is a point at infinity, and
$\gamma$ is a geodesic ray from $x$ to $z$. By a \emph{horosphere
through x with center z} we mean the limit as $t\rightarrow\infty$
of the sphere of radius $t$ in $X$ with center $\gamma(t)$. A
horosphere is the level surface of the horofunction $h(x)$
corresponding to the ray $\gamma$. By the \emph{radius} of a
horosphere through $x$ we mean the value $h(x)$. A \emph{horoball}
is the interior of a horosphere.

Suppose a group $G$ admits a properly discontinuous isometric
action on $X$ so that the quotient space $Y=X/G$ is
quasi-isometric to the union of $k$ copies of $[0,\infty)$ joined
at zero. Assume that the action of $G$ on $X$ is free. Lift the
rays in $Y$ to the rays $\gamma_i\colon[0,\infty)\rightarrow X$
for $i=1,2,\dots,k$. Let $H_i$ be the isotropy subgroup of
$\gamma_i(\infty)$; assume that $H_i$ preserves $h_i$. Assume that
in $X$ there exists a $G$-invariant system $GB$ of disjoint
horoballs, and the action of $G$ on $X\setminus GB$ is cocompact.

Then Gromov calls $G$ hyperbolic relative to the subgroups
$H_1,\dots,H_k$.

\subsection{Definitions proposed by Farb}
\begin{defn}\cite{Farb}({\bf Relatively hyperbolic group in a weak
sense}) \label{t:rhweak} Let $G$ be a finitely generated group,
and let $H$ be a finitely generated subgroup of $G.$ Fix a set $A$
of generators of $G.$ In the Cayley graph $\Gamma(G,A)$ add a
vertex $v(gH)$ for each left coset $gH$ of $H,$ and connect
$v(gH)$ with each $x\in gH$ by an edge of length $\frac 12.$ The
obtained graph $\hat{\Gamma}$ is called a coned-off graph of $G$
with respect to $H.$

The group $G$ is called \emph{weakly hyperbolic relative to H} if
$\hat{\Gamma}$ is a hyperbolic metric space.
\end{defn}

\textbf{Remark.} The terminology in Definition~\ref{t:rhweak} was
part of what was suggested by Bowditch in~\cite{Bow}. In
\cite{Farb} and \cite{Sz} a group $G$ that satisfies
Definition~\ref{t:rhweak} is termed simply ``hyperbolic relative
to H".
\begin{thm}\label{t:tSz}\cite[Theorem 1]{Sz}
Let $G$ be a finitely
generated group, and let $H_1,\dots,H_r$ be a finite set of
finitely generated subgroups of $G.$ If $G$ is hyperbolic relative
to $H_1,\dots,H_r$ in the sense of Gromov's definition, then $G$
is weakly hyperbolic relative to $H_1,\dots,H_r$.
\end{thm}
Furthermore, in \cite[Example~3]{Sz} Szczepanski shows that the
class of weakly relatively hyperbolic groups is strictly larger
than the class of groups relatively hyperbolic in the sense of
Gromov's definition.

In \cite{Farb} Farb shows solvability of the word problem for
weakly relatively hyperbolic groups that have a property which he
calls the Bounded Coset Penetration property. The Bounded Coset
Penetration property appears to be crucial for solvability of the
conjugacy problem for relatively hyperbolic groups \cite{Bu}.
\begin{defn}\label{t:BCPp}\cite{Farb}({\bf Bounded Coset Penetration
property}) Let a group $G$ be weakly hyperbolic relative to a
finitely generated subgroup $H.$ A path $u$ in $\Gamma$ is a
\emph{relative $P$-quasigeodesic}, if its projection $\hat{u}$ to
$\hat{\Gamma}$ is a $P$-quasigeodesic. The path $u$ is a
\emph{path without backtracking} if $u$ never returns to a subset
which $u$ penetrates. The pair $(G,H)$ is said to satisfy the {\it
Bounded Coset Penetration (BCP)} property if $\forall P\ge 1,$
there is a constant $c=c(P)$ so that for every pair $u,v$ of
relative $P$-quasi-geodesics without backtracking, with same
endpoints, the following conditions hold:
\begin{enumerate}
\item \label{t:bcp1} If $u$ penetrates a coset $gH$ and $v$ does
not penetrate $gH,$ then $u$ travels a $\Gamma$-distance of at
most $c$ in $gH.$ \item \label{t:bcp2} If both $u$ and $v$
penetrate a coset $gH,$ then the vertices in $\Gamma$ at which $u$
and $v$ first enter (last exit) $gH$ lie a $\Gamma$-distance of at
most $c$ from each other.
\end{enumerate}
\end{defn}
Thus one considers the class of groups which can be defined as
follows.
\begin{defn}\label{t:rhfarb}({\bf Relatively hyperbolic group by Farb})
Let $G$ be a finitely generated group, and let $H$ be a finitely
generated subgroup of $G.$ We say that $G$ is \emph{hyperbolic
relative to H in the sense of Farb}, if $G$ is weakly hyperbolic
relative to $H$ and the pair $(G,H)$ has the BCP property.
\end{defn}

 In~\cite{Sz} Szczepanski asks the following question.
\begin{quest}\label{t:qSz}\cite[Question (1)]{Sz} Let a group
$G$ be hyperbolic relative to a finitely generated subgroup $H$ in
the sense of Gromov's definition. Does the pair $(G,H)$ have the
BCP property?
\end{quest}

It is already known that Question~\ref{t:qSz} has an affirmative
answer. In fact, results of Bowditch \cite{Bow} and Dahmani
\cite{Dacl},\cite{Dadi} imply that the class of groups that are
relatively hyperbolic in the sense of Gromov's definition and the
class of groups that are relatively hyperbolic in the sense of
Farb's definition coincide.
%
%
We provide a direct proof of the affirmative answer to
Question~\ref{t:qSz}. Combined with the proof of
Theorem~\ref{t:tSz} by Szczepanski, this gives a short alternate
proof of
the following theorem.
\begin{thm}\label{t:BCPyes} Let $G$ be a finitely generated group,
and let $H$ be a finitely generated subgroup of $G.$ If the group
$G$ is hyperbolic relative to $H$ in the sense of the Gromov
definition, then $G$ is hyperbolic relative to $H$ in the sense of
the definition of Farb.
\end{thm}

\section{Quasiconvex subsets in a hyperbolic space}
Our approach is to generalize the arguments that Farb used to
prove \cite[Theorem~4.11]{Farb} (cf. \cite[Theorem~5.10]{Reb}).

\begin{notation}\label{t:not} Let $(X,d)$ be a $\delta$-hyperbolic metric space $(\delta\ge
0)$ with a collection $\Sigma$ of closed disjoint
$\epsilon$-quasiconvex subsets. We assume that the distance
between any two subsets $S_1,S_2\in\Sigma$ is bounded below by a
constant $R>24\delta+4\epsilon$. We obtain a space $X_S$ by
deleting the interiors of all of the sets in $\Sigma$. The
boundary of $X_S$ consists of disjoint connected components, each
component is the boundary of a set $S\in\Sigma$. We give $X_S$ the
path metric $d_S$. Next, we obtain the quotient $\hat{X}$ of $X_S$
by identifying points which lie in the same boundary component of
$X_S$. Hence $\hat{X}$ is equipped with a path pseudometric
$\hat{d}$ induced from the path metric $d_S$.
\end{notation}
\begin{lem}\label{t:prop1}\cite[Proposition 1]{Sz}. The space
$(\hat{X},\hat{d})$ is hyperbolic in the Gromov sense.
\end{lem}
\begin{defn}\label{defn:projection} (\textbf{Projections onto quasiconvex sets})
Fix a subset $S\in\Sigma$. Let $x$ and $y$ be two points in $X$,
and let $Pr(x)$ and $Pr(y)$ denote the projections of these points
onto $S$. Let $\alpha$ be a path in $X$ with the initial point
$i(\alpha)=x$ and the terminal point $t(\alpha)=y$. By the
\emph{projection} $Pr(\alpha,S)$ \emph{of $\alpha$ onto} $S$ we
mean an $X$-geodesic connecting $Pr(x)$ and $Pr(y)$, and by the
$X$-\emph{length $l_X(Pr(\alpha,S))$ of the projection}
$Pr(\alpha,S)$ we mean the $X$-length of that geodesic:
\[
l_X(Pr(\alpha,S))=d(Pr(x),Pr(y)).
\]
Given $S_1\in \Sigma$, by the \emph{length of the projection of
$S_1$ onto} $S$ we mean the maximum length of the projection of a
path in $S_1$ onto $S$:
\[
l_X(Pr(S_1,S))=\max\{l_X(Pr(\alpha,S))\mid \alpha\subset S_1\}.
\]
\end{defn}
\begin{notation}\label{t:projintoX}
Let $\gamma$ be a path in $X$, and let $l_X(\gamma)$ be the
$X$-length of $\gamma$. We denote by $\hat{\gamma}$ the projection
of $\gamma$ into $\hat{X}$, and by $l_{\hat{X}}(\hat{\gamma})$ the
$\hat{X}$-length of $\hat{\gamma}$. Obviously,
$l_{\hat{X}}(\hat{\gamma})\leq l_X(\gamma)$. Also, given a path
$\hat{\gamma}$ in $\hat{X}$, by $\gamma$ we mean a path in $X$
whose projection into $\hat{X}$ is $\hat{\gamma}$. Given
$\hat{\gamma}$, the path $\gamma$ is not unique in general, but
the following definition does not depend on particular choice of
$\gamma$.
\end{notation}
\begin{defn}\label{defn:intersection} (\textbf{Intersections with quasiconvex sets})
We will say that $\hat{\gamma}$ (or $\gamma$) \emph{intersects} a
subset $S\in\Sigma$, if $\gamma$ intersects the boundary of $S$
and travels a non-zero $X$-distance in the interior of $S$.
\end{defn}
By \cite[Lemma~7.3D]{Gr} and \cite[Chapter~10,
Proposition~2.1]{CDP}, we have the inequality
\begin{equation}\label{e:CDP}
d(Pr(x),Pr(y))\leq\max(C,C+d(x,y)-d(x,Pr(x))-d(y,Pr(y))),
\end{equation}
where $C=2\epsilon+12\delta$. Moreover, in~\cite{Sz} Szczepanski
shows the following. If the geodesic segment $Pr(\alpha,S)$ does
not intersect the $2\delta$-neighborhood of $\alpha$ in $X$, then
\begin{equation}\label{e:CDPSz}
l_X(Pr(\alpha,S))\leq C.
\end{equation}

Assume that each $S\in\Sigma$ is a convex set. In the proof of
Theorem~\ref{t:BCPyes} (Section~\ref{s:proof} below) we use the
fact that horoballs are convex sets. Let $\partial S$ denote the
boundary of a set $S\in \Sigma$. Observe that $\partial S$ may be
not convex (for instance, horospheres are not convex). For a set
$U\subset X$, denote by $Nb_X(U,\lambda)$ the
$\lambda$-neighborhood of $U$ in $X$.
\begin{lem}\label{lem:Sdistance} Let $S\in\Sigma$ be a set, and let
$\beta$ be an $X$-geodesic that does not intersect
$Nb_X(S,2\delta)$. If $x$ and $y$ are the endpoints of $\beta$,
then
$$d_S(Pr(x),Pr(y))\leq C+16\delta.$$
\end{lem}
\begin{proof} Let $x_s=Pr(x,S)$ and $y_s=Pr(y,S)$, so that $[x,x_s]\cap
S=x_s$ and $[y,y_s]\cap S=y_s$. By the inequality~\refe{CDPSz},
$d(x_s,y_s)\leq C$. Moreover, $X$-geodesic $[x_s,y_s]$ stays
$4\delta$-close to the union of the geodesics $[x_s,x]\cup
\beta\cup[y,y_s]$. Hence, there exists a path in $X\setminus S$
which joins $x_s$ and $y_s$, with $X$-length bounded by
$C+16\delta$, as claimed.
\end{proof}
\begin{lem}\label{lem:closeqgeod} Let $\hat{\gamma}$ be a
$P$-quasigeodesic ($P\geq 1$) in $\hat{X}$. Assume that
$\hat{\gamma}$ does not intersect any subset $S\in\Sigma$. Given a
subset $S_0\in\Sigma$, the projection $Pr(\gamma,S_0)$ of $\gamma$
onto $S_0$ has $d_S$-length at most $D=D(P)$.
\end{lem}
\begin{proof} Let $x$ and $y$ denote the endpoints of
$\hat{\gamma}$, and let $\beta$ be an $X$-geodesic joining $x$ and
$y$. Observe that $\gamma=\hat{\gamma}$, so that
$l_X(\gamma)=l_{\hat{X}}(\hat{\gamma})$. Since $\hat{d}(x,y)\leq
d(x,y)$, $\gamma$ is a $P$-quasigeodesic in $X$. If $\beta$ does
not intersect $Nb_X(S,2\delta)$, then the claim follows from
Lemma~\ref{lem:Sdistance}. In what follows, we assume that $\beta$
intersects $Nb_X(S,2\delta)$. Let $z$ (or $w$) be the point where
$\beta$ first enters (or last exits) $Nb_X(S,2\delta)$, and let
$z_s=Pr(z,S_0)$ and $w_s=Pr(w,S_0)$. It remains to show that
$d_S(z_s,w_s)$ is bounded. We argue as follows.

The Hausdorff $X$-distance between $\beta$ and $\gamma$ is bounded
by a constant $N(P)$, hence there are points $a,b\in\gamma$ so
that $d(a,z)\leq N(P)$ and $d(b,w)\leq N(P)$. Denote by $\gamma_0$
the subsegment of $\gamma$ between $a$ and $b$. Having projected
$\gamma$ to $\hat{X}$, we have the following inequalities:
\begin{align*}
\hat{d}(a,b) &\leq d(a,Pr(a,S))+d(b,Pr(b,S))\leq
2(N(P)+2\delta),\quad \text{so that}\\
l_{\hat{X}}(\hat{\gamma_0}) &\leq 2P(N(P)+2\delta).
\end{align*}
As $l_X(\gamma_0)=l_{\hat{X}}(\hat{\gamma}_0)$ and
$d_S(z_s,w_s)\leq d(z_s,a)+l_X(\gamma_0)+d(b,w_s)$, we have that
\begin{equation}\label{e:closetoS}
d_S(z_s,w_s)\leq 2(N(P)+2\delta)+2P(N(P)+2\delta).
\end{equation}
Since Lemma~\ref{lem:Sdistance} applies to the segments $[x,z]$
and $[w,y]$, we conclude that the $d_S$-length of the projection
$Pr(\gamma,S_0)$ of $\gamma$ onto $S_0$ is bounded by
\begin{equation}\label{e:D}
D(P)=2(C+16\delta)+2(P+1)(N(P)+2\delta)).
\end{equation}
\end{proof}
%
%
For each $P$-quasi-geodesic $\hat{\gamma}$ in $\hat{X}$, we are
able to bound the $d_S$-length of the projection of $\hat{\gamma}$
onto a quasiconvex set in terms of the length
$l_{\hat{X}}(\hat{\gamma})$.
\begin{lem}\label{lem:anyqgeod} Fix a subset $S_0\in\Sigma$.
Let $P\geq 1$ and $\hat{\gamma}$ be a $P$-quasigeodesic in
$\hat{X}$, which does not intersect $S_0$. Then the projection
$Pr(\gamma,S_0)$ of $\gamma$ onto $S_0$ has $d_S$-length at most
\begin{equation}\label{e:E}
l_S(Pr(\gamma),S_0)\leq
E(\hat{\gamma},P)=\big(\frac{2l_{\hat{X}}(\hat{\gamma})}{R}+3\big)(C+16\delta+D(P)).
\end{equation}
\end{lem}
\begin{proof} Let $S\in \Sigma$ be a subset different from $S_0$.
Since $S$ and $S_0$ stay distance at least $R$ apart, by
Lemma~\ref{lem:Sdistance}, the projection of $S$ onto $S_0$ has
$d_S$-length bounded by $C+16\delta$. Moreover, $\hat{\gamma}$
intersects at most $n=\frac{l_{\hat{X}}(\hat{\gamma})}{R}+1$
quasiconvex subsets from $\Sigma$, and has at most $n+2$
subsegments which do not intersect any $S\in \Sigma$. By
Lemma~\ref{lem:closeqgeod}, the claim follows.
\end{proof}
\begin{lem}\label{lem:lemma1}\cite[Lemma~1]{Sz} For any
$K\geq 2\delta+\epsilon$, there is a constant $L=L(R,K)$ so that
in the $\hat{X}$-metric, whenever $\hat{\alpha}$ is an
$\hat{X}$-geodesic, $\hat{\alpha}$ stays $(K+L/2)$-close to an
$X$-geodesic $\beta$ with the same endpoints as $\hat{\alpha}$.
\end{lem}
As a consequence of Lemma~\ref{lem:anyqgeod} and
Lemma~\ref{lem:lemma1}, we have a uniform bound on the projection
length of an $\hat{X}$-geodesic, which does not depend on the
length of this geodesic.
\begin{lem}\label{lem:fargeod} Let $K\geq
2\delta+\epsilon$, and let $\hat{\alpha}$ be a $\hat{X}$-geodesic
so that $\alpha$ does not intersect the $K$-neighborhood
$Nb_X(S_0,K)$ of $S_0$ in $X$. Then the $d_S$-length of
$Pr(\alpha,S_0)$ is bounded by a constant $E_1=E_1(R,K,\delta)$
which does not depend on the length of $\hat{\alpha}$.
\end{lem}
\begin{proof} Whenever $a\in X$ is a point, by $a_s$ we mean the projection of $a$
onto $S_0$. Let $x$ and $y$ be the endpoints of $\hat{\alpha}$,
and let $\beta$ be an $X$-geodesic joining $x$ and $y$. Observe
that if $\beta$ does not intersect $Nb_X(S_0,2\delta)$, then by
Lemma~\ref{lem:Sdistance}, $d_S(x_s,y_s)\leq C+16\delta$. In what
follows, we assume that $\beta$ intersects $Nb_X(S_0,2\delta)$.
Let $z$ (or $w$) be the point where $\beta$ first enters (or last
exits) $Nb_X(S_0,2\delta)$. By Lemma~\ref{lem:lemma1},
$\hat{\alpha}$ stays $(K+L/2)$-close to $\beta$ in the
$\hat{X}$-metric. Therefore, there are points $a,b\in\hat{\alpha}$
so that $\hat{d}(a,z)\leq K+L/2$ and $\hat{d}(b,w)\leq K+L/2$.
Hence,
\begin{equation} \label{e:alpha0}
l_{\hat{X}}(\hat{\alpha}_0)\leq 2(K+L/2+2\delta),
\end{equation}
where $\hat{\alpha}_0$ is the segment of $\hat{\alpha}$ joining
$a$ and $b$. By Lemma~\ref{lem:anyqgeod}, we have that
\begin{equation}\label{e:Salpha0}
d_S(a_s,b_s)\leq E(\hat{\alpha}_0)\leq E(2(K+L/2+2\delta),1).
\end{equation}
Let $\hat{\alpha}_1$ (or $\hat{\alpha}_2$) be the segment of
$\hat{\alpha}$ joining $x$ and $a$ (or $b$ and $w$). It remains to
show that the projections of $\hat{\alpha}_1$ and of
$\hat{\alpha}_2$ have bounded $d_S$-length. Let $\beta_i$ be an
$X$-geodesic joining the endpoints of $\hat{\alpha}_i$. As long as
$\beta_1$ intersects $Nb_X(S_0,2\delta)$, we are able to find a
point $a_1$ in $\hat{\alpha}_1$ which is $(K+L/2+2\delta)$-close
to $S_0$, so that the $\hat{X}$-length of the segment of
$\hat{\alpha}$ joining $a_1$ and $b$ satisfies the
inequality~\refe{alpha0} above. Therefore, without loss of
generality, we can assume that neither $\beta_1$ nor $\beta_2$
intersects $Nb_X(S_0,2\delta)$. Hence Lemma~\ref{lem:Sdistance}
applies to these segments, so that the length of the projection of
$\hat{\alpha}$ is bounded as follows:
\begin{equation} \label{e:pralpha}
l_S(Pr(\alpha),S_0) \leq E_1=E(2K+L+4\delta,1)+2(C+16\delta).
\end{equation}
This upper bound does not depend on the length of $\hat{\alpha}$,
as claimed.
\end{proof}
\begin{lem}\label{t:central} \textbf{(Bounded Subset Penetration)}
Let $P\geq 1$ be a constant, and let $\hat{\xi}$ and $\hat{\tau}$
be two $P$-quasigeodesics without backtracking in $\hat{X}$, with
common endpoints. Then there exists a constant
$B=B(\delta,\epsilon,P)$ such that the following holds.
\begin{enumerate}
  \item\label{e:first} If $\hat{\xi}$ intersects a subset
  $S\in\Sigma$ and $\hat{\tau}$ does not intersect $S$, then the $d_S$-distance
  $s$ between the points where $\xi$ first enters and last exits
  $S$ is bounded by $B$.
  \item\label{e:second} If both $\hat{\xi}$ and $\hat{\tau}$ intersect a subset
  $S\in\Sigma$, then the $d_S$-distance $s$ between the points where $\hat{\xi}$ and
  $\hat{\tau}$ first enter (or last exit) $S$ is bounded by $B$.
\end{enumerate}
\end{lem}
\begin{proof} Fix $K=2\delta+\epsilon$. In the
case~(\ref{e:first}), observe that $s$ is bounded by the
$d_S$-length of the projection of the path which is the
concatenation of $\hat{\tau}$ and the two subsegments of
$\hat{\xi}$ that lie outside $S$. Let $\hat{\alpha}$ denote either
of these $P$-quasigeodesic segments, and let $\hat{\beta}$ be a
$\hat{X}$-geodesic with the same endpoints as $\hat{\alpha}$.
Observe that if a subsegment $\hat{\beta}'$ of $\hat{\beta}$ does
not intersect $Nb_{\hat{X}}(S,K)$, then $\beta'$ does not
intersect $Nb_X(S,K)$, hence by Lemma~\ref{lem:fargeod},  we have
that $l_S(Pr(\beta'))\leq E_1$.

Next , assume that $\hat{\beta}$ intersects $Nb_X(S,K)$, but does
not intersect $S$. Let $x$ be the point where $\hat{\beta}$ first
enters $Nb_{\hat{X}}(K,S)$, and let $y$ be the point where
$\hat{\beta}$ last exits $Nb_{\hat{X}}(S,K)$. Denote by
$\hat{\beta}_{xy}$ the segment of $\hat{\beta}$ joining $x$ and
$y$. Then $l_{\hat{X}}(\hat{\beta}_{xy})\leq \hat{d}(x,y)\le 2K$.
Since $R>3K$, the geodesic $\hat{\beta}_{xy}$ does not leave
$Nb_{\hat{X}}(S,R)$, so that $\hat{\beta}_{xy}$ does not intersect
any subset $S'\in\Sigma$, and by Lemma~\ref{lem:closeqgeod},
$l_S(Pr(\beta_{xy},S))\leq D(1)$. Observe that $\hat{\beta}$ is
the concatenation of $\hat{\beta}_{xy}$ and other two (possibly,
degenerate) segments that do not intersect $Nb_{\hat{X}}(S,K)$.
Therefore, we have that $l_S(Pr(\hat{\beta},S))\leq D(1)+2E_1$ in
this case.

Finally, assume that $\hat{\beta}$ intersects $S$. Let $x_s$ (or
$y_s$) be the point where $\hat{\beta}$ first enters (or last
exits) $S$. As $\hat{X}$ is a hyperbolic space, $\hat{\alpha}$ and
$\hat{\beta}$ stay a bounded distance apart; let $M(P)$ be a
constant that bounds this distance. There are points
$z,w\in\hat{\alpha}$ with $\hat{d}(z,x_s)\leq M(P)$ and
$\hat{d}(w,y_s)\leq M(P)$. Let $\hat{\alpha}_{zw}$ be the segment
of $\hat{\alpha}$ between $z$ and $w$, we have that
$l_{\hat{X}}(\hat{\alpha}_{zw})\leq 2PM(P)$. By
Lemma~\ref{lem:anyqgeod},
\[
l_S(Pr(\alpha_{zw}),S)\leq E(\hat{\alpha}_{zw})\leq E(2PM(P)).
\]
Let $\hat{\alpha}_1$ (or $\hat{\alpha}_2$) be the segment of
$\hat{\alpha}$ joining the initial point of $\hat{\alpha}$ and $z$
(or $w$ and the terminal point of $\hat{\alpha}$). Let $\beta_i$
be an $X$-geodesic joining the endpoints of $\hat{\alpha}_i$.
W.l.o.g. (cf. the proof of Lemma~\ref{lem:fargeod}), we can assume
that $\beta_1,\beta_2$ do not intersect $S$. As we have shown,
$l_S(Pr(\beta_i,S))\leq D(1)+2E_1$, for $i=1,2$.

Therefore, $l_S(Pr(\alpha,S))\leq 2(D(1)+2E_1)+E(2PM(P))$, so that
\begin{equation}\label{e:BCP1}
s\leq 3(2(D(1)+2E_1)+E(2PM(P)))
\end{equation}
which finishes the proof in the case~\refe{first}.

Now, we prove~\refe{second}. Let $s$ be the $X$-distance between
the points where $\xi$ and $\tau$ first enter the subset $S$, and
let $\xi_1$ and $\tau_1$ be the initial segments of $\xi$ and
$\tau$ with ends at the points where these paths first enter $S$.
Observe that $s\leq l_S(Pr(\xi_1^{-1}\circ\tau_1),S)$. The
arguments used in the proof of case~\refe{first} show that
\[
s\leq 2(2(D(1)+2E_1)+E(2PM(P))).
\]
Obviously, the distance between the points where $\xi$ and $\tau$
last exit the subset $S$ is bounded by $2(2(D(1)+2E_1)+E(2PM(P)))$
as well.

Therefore, one can set $B=3(2(D(1)+2E_1)+E(2PM(P)))$ to finish the
proof of the lemma.
\end{proof}

\section{Proof of Theorem~\ref{t:BCPyes}}\label{s:proof}
First, consider the case of one subgroup $H$ so that $G$ is
hyperbolic relative to $H$. By \cite[Corollary 3.2 ]{Farb}, we can
assume that the generating set $A$ of $G$ contains a generating
set $A_H$ of $H$ as a subset. We define a map $f$ from the Cayley
graph $\Gamma(G,A)$ to $X$ as follows. Lift the ray in the
quotient space $Y$ to a ray $\gamma\colon [0,\infty)\rightarrow X$
and choose a horoball $S$ so that the boundary of $S$ is the level
surface of the horofunction which corresponds to the ray $\gamma$,
and the images $GS$ of $S$ under the action of $G$ form a
$G$-invariant system of disjoint horoballs. We assume that the
minimum distance between two horoballs is $R$. Notifying that a
horoball is a convex subset of $X$, we define $X'$ to be the space
obtained from $X$ by deleting the interiors of all the horoballs
in $\Sigma=GS$ (cf. notation~\ref{t:not}). Clearly, the action of
$G$ on $X'$ is cocompact.

Pick a point $x$ on the boundary of $S$. Let $g.x$ be the image of
$x$ under the action of $g\in A$. Hence, $g.x\in g.S$. Observe
that the horoball $g.S$ corresponds to the left coset $gH$ in the
following meaning. If $g\in A\setminus A_H$, then $g.S\neq S$, and
we join $x$ and $g.x$ by a segment which does not intersect any
horoball. If $g\in A_H$, then $g.S=S$, and we join $x$ and $g.x$
by a segment which lies in $S$. Define a map $f\colon
\Gamma(G,A)\rightarrow X$ as follows. Let $id\in \Gamma(G,A)$ be
the point that corresponds to $1_G$. The map $f$ sends $id\in
\Gamma(G,A)$ to the point $x\in X$ and $[id,g.id]$ to the segment
in $X$ that joins $x$ and $g.x$, for each $g\in A$. Extend $f$ to
$\Gamma(G,A)$ equivariantly. Since $G$ acts on $X$ by isometries,
and the action of $G$ on $\hat{X}$ is cocompact, the map $f$
induces a quasi-isometry $\hat{f}$ between $\hat{\Gamma}$ and
$\hat{X}$. Therefore, \cite[Proposition~1]{Sz} (cf.
Lemma~\ref{t:prop1}) implies that $\hat{\Gamma}$ is a hyperbolic
metric space i.e., $G$ is weakly hyperbolic relative to $H$
(\cite[Theorem~1]{Sz}, cf. Theorem~\ref{t:tSz}).

Furthermore, let $\rho$ be a path in $\Gamma(G,A)$, and let
$f(\rho)=\beta\in X$. If $\rho$ penetrates a coset $gH$, then
$\beta$ intersects the horosphere $g.S$; moreover, if the distance
between the points where $\beta$ first enters and last exits the
boundary of $g.S$ is bounded by $l$, then the $\Gamma$-distance
that $\rho$ travels in $gH$ is bounded in terms of $l$ and the
constants of quasi-isometry $\hat{f}$. Therefore, the pair $(G,H)$
has the BCP property.

Now, assume that $G$ is hyperbolic relative to $H_1,\dots,H_k$, in
the sense of the Gromov definition. We assume that the generating
set $A$ of $G$ contains a generating set $A_i\subset A$ for $H_i$,
for all $i=1,\dots,k$. Lift the $k$ rays in the quotient space $Y$
to rays $\gamma_i\colon [0,\infty)\rightarrow X$ and choose
horoballs $S_i$ so that the boundary of $S_i$ is the level surface
of the horofunction which corresponds to the ray $\gamma_i$, and
the images $\Sigma=\bigcup_{i=1}^k GS_i$ of $S_i$ under the action
of $G$ form a $G$-invariant system of disjoint horoballs. Assume
that the distance between two horoballs in $\Sigma$ is bounded
below by $R$. Pick a point $x_i$ in $S_i$, for each $i=1,\dots,k$,
and connect each $x_j$ to $x_1$ by an edge $e_j$, for
$j=2,\dots,k$. The union of the edges $e_j$ is a finite tree which
we denote by $T$. Each generator $h\in A_i$ (for some $i$) maps
$x_i\in S_i$ to another point $y_i\in S_i$, and $x_j\in S_j$
($j\neq i$) to a point $h.x_i\in h.S_i\neq S_i$, so that $hT$ and
$T$ are disjoint. Obviously, $T$ and $g.T$ are disjoint for each
generator $g\in A\setminus\bigcup_{i=1}^k A_i$. We proceed as in
the case of one subgroup and get a map from $\hat{\Gamma}$ to
$\hat{X}$. Contract the tree $T$ and all its images $GT$ to a
point, in order to see that $\hat{f}$ is a quasi-isometry between
$\hat{\Gamma}$ and $\hat{X}$. The same argument as above shows
that the BCP property holds.

\section*{Acknowledgements} This paper was written during my stay at
McGill University as a postdoctoral fellow. I am thankful to my
supervisors O.~Kharlampovich, A.~Miasnikov and D.~Wise for their
warm hospitality. I am also thankful to Fran\c{c}ois Dahmani for
kindly submitting me his thesis, and to the referee for useful
comments.

\end{document}